\documentclass[11pt]{amsart}
\usepackage[a4paper]{geometry}
\usepackage{color}
\usepackage{longtable}

\usepackage{lipsum}
\newcommand{\mylongtitle}[1]{%
  \ifodd\value{page}%
    \protect\parbox{0.9\linewidth}{#1}\hfill%
  \else%
    \hfill\protect\parbox{0.9\linewidth}{#1}%
  \fi%
}

 \newtheorem{thm}{Theorem}[section]
 
 \newtheorem{lem}[thm]{Lemma}
 
 \theoremstyle{definition}
 
 \theoremstyle{remark}
 
 \numberwithin{equation}{section}

\begin{document}
\title {\textbf{ \mylongtitle{Lorentz Hypersurfaces satisfying $\triangle \vec {H}=  \alpha \vec {H}$ with non diagonal shape operator}}}
\author{Deepika, Andreas Arvanitoyeorgos and Ram Shankar Gupta}

\begin{abstract}
We study Lorentz hypersurfaces $M_{1}^{n}$ in
$E_{1}^{n+1}$ satisfying $\triangle \vec {H}=  \alpha \vec {H}$ with
non diagonal shape operator, having complex eigenvalues. We prove
that every such Lorentz hypersurface in $E_{1}^{n+1}$ having at most
five distinct principal curvatures has constant mean curvature.
\\
\\
\textit{AMS 2000 MSC Classification.} Primary 53D12 Secondary 53C40; 53C42
\end{abstract}

\maketitle
\section{\textbf{Introduction}}

The study of submanifolds with harmonic mean curvature vector field
was initiated by B.Y. Chen in 1985 and arose in the context of his
theory of submanifolds of finite type. For a survey on submanifolds
of finite type and various related topics we refer to \cite
{BY2,BY3}. Let $M^{n}_{r}$ be an $n$-dimensional, connected
submanifold of the pseudo-Euclidean space $E^{m}_{s}$. We denote by
$\vec {x}$, $\vec {H}$, and $\triangle $ respectively the position
vector field, mean curvature vector field of $M^{n}_{r}$, and the
Laplace operator on $M^{n}_{r}$, with respect to the induced metric
$g$ on $M^{n}_{r}$, from the indefinite metric on the ambient space
$E^{m}_{s}$. It is well known (\cite {BY1}) that

\begin{equation} \label{1.1}
\triangle \vec {x} = -n \vec {H}.
\end{equation}

 A submanifold $M^{n}_{r}$ of $E^{m}_{s}$ satisfying the condition

\begin{equation} \label{1.2}
\triangle \vec {H} = 0,
\end{equation}
is called {\it biharmonic} submanifold. In view of (\ref{1.1}),
condition (\ref {1.2}) is equivalent to $\triangle ^{2}\vec {x}=0$.
Equation (\ref {1.2}) is a special case of the equation

\begin{equation}\label{1.3}
\triangle \vec {H}=  \alpha \vec {H}.
\end{equation}

 As remarked above, minimal submanifolds are immediately seen to be
biharmonic. Conversely, a question arises whether the class of
submanifolds with harmonic mean curvature vector field is
essentially larger than the class of minimal submanifolds.
Concerning this problem B.Y. Chen conjectured
the following:\\

\noindent
 \emph{\textbf{Conjecture.} The only biharmonic submanifolds of Euclidean
spaces are the minimal ones.}\\

In Euclidean spaces, we have the following results, which indeed
support the above mentioned conjecture. B.Y. Chen proved in 1985
that every biharmonic surface in $E^{3}$ is minimal. Thereafter, I.
Dimitric generalized this result in \cite {ID}. In \cite {TT}, it
was proved by Th. Hasanis and Th. Vlachos that every biharmonic
hypersurface in $E^{4}$ is minimal. Recently, it was proved by the
third author that every biharmonic hypersurface with three distinct
principal curvatures in $E^{n+1}$ with arbitrary dimension is
minimal (\cite {RS}).

The study of equation (\ref {1.3}) for submanifolds in
pseudo-Euclidean spaces was originated by Ferrandez et al. in \cite
{AP1,AP2}. They showed that if the minimal polynomial of the shape
operator of a hypersurface $M^{n-1}_{r}$  $(r = 0, 1)$ in
$E^{n}_{1}$ has degree at most two, then $M^{n-1}_{r}$ has constant
mean curvature. Also, in \cite {BY2} various classification theorems
for submanifolds in a Minkowski spacetime were obtained. In \cite
{AFG}, it was proved that every hypersurface $M^{3}_{r}$ $(r = 0, 1,
2, 3)$ of $E^{4}_{s}$ satisfying equation (\ref {1.3}) whose shape
operator is diagonal, has constant mean curvature. Also, in \cite
{AGM} the same conclusion was obtained for every hypersurface
$M^{3}_{1}$ in $E^{4}_{1}$. In \cite {JC}, it was proved that every
hypersurface having at most three distinct principal curvatures in
$E^{n+1}_{s}$ satisfying (\ref {1.3}) with diagonal shape operator
has constant mean curvature. Recently, it was proved that every
biharmonic hypersurfaces in $E^{5}$ with all distinct principal
curvatures is minimal \cite {RSS}.

In contrast to the submanifolds of Euclidean spaces, Chen's
conjecture is not always true for the submanifolds of the
pseudo-Euclidean spaces. For example, B.Y. Chen et al. \cite
{BS1,BS2} obtained some examples of proper biharmonic surfaces in
4-dimensional pseudo-Euclidean spaces $E^{4}_{s}$  for $s = 1, 2, 3$
(see also \cite {BY4}). However, it is reasonable to expect that for
hypersurfaces in pseudo-Euclidean space, Chen's conjecture is true.
This is supported by the following facts: B. Y. Chen et al. proved
in \cite {BS1,BS2} that biharmonic surfaces in pseudo-Euclidean
3-spaces are minimal, and the second author et al. \cite {AFGV}
proved that biharmonic Lorentzian hypersurfaces in Minkowski
4-spaces are minimal. Recently, it was proved that every Lorentz
hypersurface $M_{1}^{n}$ in $E^{n+1}_{1}$ having complex eigenvalues
with at most four distinct principal curvatures has constant mean
curvature \cite {DR}.

  In this paper, we study Lorentz hypersurfaces
$M_{1}^{n}$ in $E_{1}^{n+1}$  satisfying (\ref {1.3}) and having
shape operator (\ref {2.11}) with at most five distinct eigenvalues.
Our main result is the following:

\begin{thm}\label{main}
Every Lorentz hypersurface $M_{1}^{n}$ in $E^{n+1}_{1}$ satisfying
$\triangle \vec {H}=  \alpha \vec {H}$, having non diagonal shape
operator with complex eigenvalues with at most five distinct
principal curvatures, has
constant mean curvature.
\end{thm}

\noindent
We briefly present the central ideas of the
proof of the above theorem

\noindent We use the proper mean curvature condition and successive
use of the Codazzi equation to simplify the connection forms
defining the covariant derivative.  In this way  we obtain relations
among the eigenvalues of the shape operator, the connection forms
and the mean curvature $H$ (cf. Table 1). Next, we use the Gauss
equation and covariant differentiation with respect to an
orthonormal frame to prove that the real part of complex eigenvalues
vanishes, and obtain that the mean curvature $H$ is either zero or
constant.

This is the most difficult part of the proof and it
  is achieved by showing that $H$ satisfies a non trivial algebraic polynomial equation
with constant coefficients, hence it must be constant.
Reaching to such a polynomial equation is not a trivial matter in general, and unfortunately it seems there is no standard method to get it.
In our case we obtain two polynomials with coefficients in the polynomial ring $R[H]$ that have one of the eigenvalues of the shape operator as a root.  Then, by using standard argument involving the resultant of two polynomials, it follows that  $H$ must be constant.

\section{\textbf{Preliminaries}}

Let ($M_{1}^{n}, g$) be a $n$-dimensional Lorentz hypersurface
isometrically immersed in a $n+1$-dimensional pseudo-Euclidean space
($E^{n+1}_{1}, \overline g$) and $g = \overline g_{|M_{1}^{n}}$. We
denote by $\xi$ the unit normal vector to $M_{1}^{n}$ with
$\overline g(\xi, \xi)= 1$.

  Let $\overline\nabla $ and $ \nabla $ denote the linear connections on $E_{1}^{n+1}$ and $M_{1}^{n}$ respectively. Then, the Gauss and Weingarten formulae are given by
\begin{equation}\label{2.1}
\overline\nabla_{X}Y = \nabla_{X}Y + h(X, Y), \hspace{.3 cm} X, Y
\in\Gamma(TM^{n}_{1}),
\end{equation}
\begin{equation}\label{2.2}
\overline\nabla_{X}\xi = -S_{\xi}X, \hspace{.3 cm}  \xi
\in\Gamma(TM^{n}_{1})^{\bot},
\end{equation}
where $h$ is the second fundamental form and $S$ is the shape
operator. It is well known that $h$ and
$S$ are related by

\begin{equation}\label{2.3}
\overline{g}(h(X,Y), \xi) = g(S_{\xi}X,Y).
\end{equation}

The mean curvature vector is given by
\begin{equation}\label{2.4}
\vec{H} = \frac{1}{n} \rm trace \emph{h}.
\end{equation}

The Gauss and Codazzi equations are given by
\begin{equation}\label{2.5}
R(X, Y)Z = g(SY, Z) SX - g(SX, Z) SY,
\end{equation}

\begin{equation}\label{2.6}
(\nabla_{X}S)Y = (\nabla_{Y}S)X
\end{equation}
respectively, where $R$ is the curvature tensor, $S=S_{\xi}$ for
some unit normal vector field $\xi$ and
\begin{equation}\label{2.7}
(\nabla_{X}S)Y = \nabla_{X}(SY)- S(\nabla_{X}Y),
\end{equation}
for all $ X, Y, Z \in \Gamma(TM_{1}^{n})$.

The necessary and sufficient
 conditions for $M^{n}_{1}$ to have proper mean curvature in $E^{n+1}_{1}$
 are (\cite {AFG})
\begin{equation}\label{2.8}
\triangle H + H \rm trace \emph{S}^{\emph{2}} = \alpha \emph{H},
\end{equation}
\begin{equation}\label{2.9}
 S (\rm grad \emph{H})+ \frac{n}{2} \emph{H} grad \emph{H} = 0,
\end{equation}
where $H$ denotes the mean curvature. Also, the Laplace operator
$\triangle$ of a scalar valued function $f$ is given by (\cite
{BS1})
\begin{equation}\label{2.10}
\triangle f = -\sum_{i=1}^{n}\epsilon_{i}(e_{i}e_{i}f -
\nabla_{e_{i}}e_{i}f),
\end{equation}
where $\{e_{1}, e_{2},..., e_{n}\}$ is an orthonormal local tangent
frame on $M^{n}_{1}$ with $\epsilon_{i}=\pm 1$.\vspace{.4cm}

   A vector $X$ in $E_{s}^{n+1}$ is called
  spacelike, timelike or lightlike according if $\overline g(X, X)> 0$,
  $\overline g(X, X)<0$ or $\overline g(X,
  X)=0$ respectively. A non degenerate hypersurface $M^{n}_{r}$ of $E_{s}^{n+1}$ is
  called Riemannian or pseudo-Riemannian according as the induced metric on
  $M_{r}^{n+1}$ from the indefinite metric on $E_{s}^{n+1}$ is definite
  or indefinite. The shape operator of
  pseudo-Riemannian hypersurfaces is not always diagonalizable in contrast to the
  Riemannian hypersurfaces.

 The matrix
representation of the shape operator of $M^{n}_{1}$ in $E_{1}^{n+1}$
having complex eigen values with respect to a suitable orthonormal
base field of the tangent bundle takes the form (\cite {AZ,MA})

\begin{equation}\label{2.11}
 S = \left(
                            \begin{array}{ccccc}
                               \lambda & -\mu & & \\
                              \mu & \lambda & \\
                              & &     D_{n-2} \\
                             \end{array}
                          \right), \hspace{.5 cm}
 \end{equation}
where $\mu \neq 0$ and
$D_{n-2}=$ ${\rm diag}(\lambda_{3},\dots,\lambda_{n})$.

 The following algebraic lemma is useful to get our result:

\begin{lem}\label{lemma2.1}\cite[Theorem 4.4, pp. 58--59]{KK} Let D be a unique factorization domain, and let $f(X) =
a_{0}X^{m} +a_{1}X^{m-1} + \cdots + a_{m}, g(X) = b_{0}X^{n} +
b_{1}X^{n-1} + \cdots + b_{n}$ be two polynomials in $D[X]$. Assume
that the leading coefficients $a_{0}$ and $b_{0}$ of $f(X)$ and
$g(X)$ are not both zero. Then $f(X)$ and $g(X)$ have a non constant
common factor if and only if the resultant $\Re(f, g)$ of $f$ and
$g$ is zero, where
\begin{center}
$\Re(f,g)=\det
\begin{pmatrix}
  a_{0} & a_{1} & a_{2} & \cdots & a_{m} &   &   &   \\
    & a_{0} & a_{1} & \cdots & \cdots & a_{m} &   &   \\
    &   & \ddots & \ddots & \ddots & \ddots & \ddots &   \\
   &   &   & a_{0} & a_{1} & a_{2} & \cdots & a_{m} \\
  b_{0} & b_{1} & b_{2} & \cdots & b_{n} &   &   &   \\
    & b_{0} & b_{1} & \cdots & \cdots & b_{n} &  &   \\
    &   & \ddots & \ddots & \ddots & \ddots & \ddots &   \\
    &   &   & b_{0} & b_{1} & b_{2} & \cdots & b_{n} \\
\end{pmatrix}.
$
\end{center}
In the above determinant there are $n$ rows of ``$a$" entries and $m$ rows of ``$b$" entries.

\end{lem}

 \section{\textbf{Lorentz hypersurfaces in $E_{1}^{n+1}$ satisfying $\triangle \vec {H}=  \alpha \vec {H}$}}

We assume that $H$ is not constant and grad$H\neq0$. Then there
exists an open connected subset $U$ of $M^{n}_{1}$, with grad$_{p}H
\neq 0$ for all $p\in U$. From (\ref {2.9}), it is easy to see that
grad$H$ is an eigenvector of the shape operator $S$ with the
corresponding principal curvature $-\frac{n}{2}H$. In view of (\ref
{2.11}), the shape operator $S$ of hypersurfaces will take the
following form
\begin{multline}\label{3.1}
S(e_{1})=\lambda e_{1}+\mu e_{2},\hspace{.2 cm} S(e_{2})=-\mu
e_{1}+\lambda e_{2},\hspace{.2 cm} S(e_{3})=\lambda_{3} e_{3},\dots,
S(e_{n})=\lambda_{n} e_{n},
 \end{multline}
 with respect to orthonormal basis $\{e_{1}, e_{2},..., e_{n}\}$ of
$T_{p} M_{1}^{n}$, which satisfies
\begin{equation}\label{3.2}
 g(e_{1},e_{1})=-1, \quad g(e_{i},e_{i})=1, \hspace{.2 cm}
 i=2,3,...,n,
 \end{equation}
and
 \begin{equation}\label{3.3}
  g(e_{i},e_{j})=0, \hspace{.2 cm} for \hspace{.2 cm}  i \neq j.
\end{equation}

We write
\begin{equation}\label{3.4}
    \nabla_{e_{i}}e_{j}=\sum_{k=1}^{n}\omega_{ij}^{k}e_{k},\hspace{2
cm} i, j = 1, 2, ... , n.
\end{equation}

Using (3.4) and taking covariant derivatives of (3.2) and (3.3) with
$e_{k}$, we find
     \begin{equation}\label{3.5}
   \begin{split}
    \omega_{ki}^{i}=0,
    \hspace{1 cm}
    \omega_{kj}^{i}=-\omega_{ki}^{j},
    \end{split}
    \end{equation}
for $i \neq j$ and  $i, j, k =1, 2, ... , n.$\\

 In view of (\ref {3.1}), grad$H$ can be chosen in one of
the directions $e_{3},\dots, e_{n}$ and in each direction grad$H$ is
space-like. Without loss of generality, we can choose $e_{n}$ in the
direction of grad$H$, so $\lambda_{n}= -\frac{nH}{2}.$  We express
grad$H$ as grad$H$ = $-e_{1}(H)e_{1}+e_{2}(H)e_{2}+ \dots
+e_{n}(H)e_{n}$, which gives

\begin{equation} \label{3.6}
e_{n}(H)\neq 0,\quad e_{1}(H)=e_{2}(H)= \cdots = e_{n-1}(H)=0.
\end{equation}

Using (\ref {3.4}), (\ref {3.6}) and the fact that $[e_{i}
\hspace{.1 cm}
e_{j}](H)=0=\nabla_{e_{i}}e_{j}(H)-\nabla_{e_{j}}e_{i}(H),$ for
$i\neq j$ and $i, j \neq n$, we find
\begin{equation}\label{3.7}
\omega_{ij}^{n}=\omega_{ji}^{n}.
\end{equation}

From (\ref {2.7}), (\ref {3.1}), (\ref {3.4}) and (\ref {3.6}), the
Codazzi equation $g((\nabla_{e_{n}}S)e_{a},e_{a}) =
g((\nabla_{e_{a}}S)e_{n},e_{a})$ leads to

\begin{equation}\label{3.8}
e_{n}(\lambda_{a})= (\lambda_{n}-\lambda_{a})\omega_{na}^{a}, \quad
3\leq a\leq n-1
\end{equation}

Therefore, $\lambda_{n}\neq \lambda_{a}$, because if $\lambda_{n}=
\lambda_{a}$, from (\ref {3.8}) we have $e_{n}(H)=0$, which
contradicts (\ref {3.6}).\\

Due to the main result of \cite[Theorem 3.2]{DR} it suffices to
consider
only the case of precisely five distinct principal curvatures.\\

\emph{From now on we assume that the shape operator (\ref {2.11})
has five distinct
eigenvalues.}\\

 It can be easily seen that the eigenvalues of the shape operator (\ref {2.11}) are $\lambda \pm \sqrt{-1} \mu,
 \lambda_{3}, \dots ,\lambda_{n}$. So, under the assumption that the shape operator
 (\ref {2.11}) has five distinct eigenvalues, we can assume that $\lambda_{3}=\lambda_{4}= \cdots =\lambda_{r}$ and $\lambda_{r+1}=\lambda_{r+2}= \cdots
 =\lambda_{n-1}$.
 So, expressions (\ref {3.1}) reduce to
\begin{equation}\label{3.9}
S(e_{1})=\lambda e_{1}+\mu e_{2}, \ S(e_{2})=-\mu
e_{1}+\lambda e_{2},\  S(e_{A})=\lambda_{3}
e_{A},\  S(e_{B})=\lambda_{n-1} e_{B},\
S(e_{n})=\lambda_{n} e_{n},
 \end{equation}
where $A=3,4, \dots ,r$ and $B= r+1, r+2, \dots ,n-1$.

From now on we assume that
   \begin{center}
   $\begin{array}{lcl} A\neq \widetilde{A}, \quad A, \widetilde{A}= 3, 4, \dots, r, \\  B\neq \widetilde{B},
  \quad B,  \widetilde{B}= r+1, r+2,\dots, n-1, .\end{array}$
   \end{center}

From $g((\nabla_{X}S)Y,Z) = g((\nabla_{Y}S)X,Z)$, using (\ref
{2.7}), (\ref {3.4}), (\ref {3.6}), (\ref{3.9}) and the value
$\lambda_{n}=-\frac{nH}{2}$, we obtain the equations Ti in Table 1
showing the relations among the eigenvalues of $S$, the connection
forms $\omega _{ij}^k$, and the  orthonormal frame $\{e_i\}$.

\smallskip
 \begin{center}
  Table 1. \ {Evaluation of $g((\nabla _XS)Y, Z)=g((\nabla _YS)X, Z)$ for various values of $e_i$.}
\end{center}

\begin{longtable}{|c | c |c |c |c |}
\hline 
i & X & Y  & Z & Codazzi equation  Ti  \\ [0.5ex] 
\hline\hline 
1 &$e_{1}$ & $e_{2}$ & $e_{1}$ & $e_{2}(\lambda) + e_{1}(\mu)=0 $\\

2 &$e_{1}$ & $e_{2}$ & $e_{2}$ & $e_{1}(\lambda) - e_{2}(\mu)=0 $\\

 3 & $e_{1}$ & $e_{2}$ & $e_{A}$ &
$[\lambda-\lambda_{3}](\omega_{12}^{A}-\omega_{21}^{A})=
\mu(\omega_{22}^{A}+\omega_{11}^{A})$ \\ 

 4 & $e_{1}$ & $e_{2}$ & $e_{B}$ & $[\lambda -
\lambda_{n-1}](\omega_{12}^{B}-\omega_{21}^{B})=
\mu(\omega_{22}^{B}+\omega_{11}^{B})$\\

 5 & $e_{1}$ & $e_{2}$ & $e_{n}$ & $[\lambda +
\frac{nH}{2}](\omega_{12}^{n}-\omega_{21}^{n})=
\mu(\omega_{22}^{n}+\omega_{11}^{n})$\\

\hline
6 & $e_{1}$ & $e_{A}$ & $e_{1}$  & $e_{A}(\lambda)= [\lambda_{3}- \lambda]\omega_{1A}^{1} + \mu \omega_{1A}^{2}$\\

 7 & $e_{1}$ & $e_{A}$ & $e_{2}$ & $e_{A}(\mu)= [\lambda_{3}-
\lambda]\omega_{1A}^{2} - \mu \omega_{1A}^{1}
$ \\

 8 & $e_{1}$ & $e_{A}$ & $e_{A}$ & $ e_{1}(\lambda_{3})= [\lambda-
\lambda_{3}]\omega_{A1}^{A} + \mu \omega_{A2}^{A}
$ \\

 9 & $e_{1}$ & $e_{A}$ & $e_{\widetilde{A}}$ & $[\lambda-
\lambda_{3}]\omega_{A1}^{\widetilde{A}} + \mu
\omega_{A2}^{\widetilde{A}}=0$\\

 10 & $e_{1}$ & $e_{A}$ & $e_{B}$ & $[\lambda_{3}-
\lambda_{n-1}]\omega_{1A}^{B} = [\lambda -
\lambda_{n-1}]\omega_{A1}^{B}
+ \mu \omega_{A2}^{B}$ \\

 11 & $e_{1}$ & $e_{A}$ & $e_{n}$ & $[\lambda_{3} +
\frac{nH}{2}]\omega_{1A}^{n} = [\lambda +
\frac{nH}{2}]\omega_{A1}^{n}
+ \mu \omega_{A2}^{n}$\\

\hline
12 & $e_{1}$ & $e_{B}$ & $e_{1}$ & $ e_{B}(\lambda)= [\lambda_{n-1}- \lambda]\omega_{1B}^{1} + \mu \omega_{1B}^{2}$\\

 13 & $e_{1}$ & $e_{B}$ & $e_{2}$ & $e_{B}(\mu)= [\lambda_{n-1}-
\lambda]\omega_{1B}^{2} - \mu
\omega_{1B}^{1}$\\

 14 & $e_{1}$ & $e_{B}$ & $e_{A}$ & $[\lambda_{n-1}-
\lambda_{3}]\omega_{1B}^{A} = [\lambda - \lambda_{3}]\omega_{B1}^{A}
+ \mu \omega_{B2}^{A}$\\

 15 & $e_{1}$ & $e_{B}$ & $e_{B}$ & $ e_{1}(\lambda_{n-1})=
[\lambda- \lambda_{n-1}]\omega_{B1}^{B} + \mu
\omega_{B2}^{B}$\\

 16 & $e_{1}$ & $e_{B}$ & $e_{\widetilde{B}}$ & $[\lambda-
\lambda_{n-1}]\omega_{B1}^{\widetilde{B}} + \mu
\omega_{B2}^{\widetilde{B}}=0$\\

 17 & $e_{1}$ & $e_{B}$ & $e_{n}$ & $[\lambda_{n-1} +
\frac{nH}{2}]\omega_{1B}^{n} = [\lambda +
\frac{nH}{2}]\omega_{B1}^{n}
+ \mu \omega_{B2}^{n}$\\

\hline
18 & $e_{1}$ & $e_{n}$ & $e_{1}$ & $ -(\lambda + \frac{nH}{2})\omega_{1n}^{1}+ \mu \omega_{1n}^{2}=e_{n}(\lambda)$\\

 19 & $e_{1}$ & $e_{n}$ & $e_{2}$ & $-(\lambda +
\frac{nH}{2})\omega_{1n}^{2}- \mu
\omega_{1n}^{1}=e_{n}(\mu)$\\

20 & $e_{1}$ & $e_{n}$ & $e_{n}$ & $ (\lambda + \frac{nH}{2})\omega_{n1}^{n}+ \mu \omega_{n2}^{n}=0$\\

\hline
21 & $e_{2}$ & $e_{A}$ & $e_{1}$ & $ -e_{A}(\mu)= [\lambda_{3}-\lambda]\omega_{2A}^{1}+ \mu \omega^{2}_{2A}$\\

22 & $e_{2}$ & $e_{A}$ & $e_{2}$ & $ e_{A}(\lambda)=[\lambda_{3}-\lambda]\omega_{2A}^{2}- \mu \omega^{1}_{2A}$\\

23 & $e_{2}$ & $e_{A}$ & $e_{A}$ & $ e_{2}(\lambda_{3})= [\lambda-\lambda_{3}]\omega_{A2}^{A}- \mu \omega^{A}_{A1}$\\

24 & $e_{2}$ & $e_{A}$ & $e_{\widetilde{A}}$ & $[\lambda-\lambda_{3}]\omega_{A2}^{\widetilde{A}}- \mu \omega^{\widetilde{A}}_{A1}=0$\\

25 & $e_{2}$ & $e_{A}$ & $e_{B}$ & $(\lambda_{3}-\lambda_{n-1})\omega_{2A}^{B}= [\lambda-\lambda_{n-1}]\omega_{A2}^{B}- \mu \omega^{B}_{A1}$\\

 26 & $e_{2}$ & $e_{A}$ & $e_{n}$ &
$[\lambda_{3}+\frac{nH}{2}]\omega_{2A}^{n} = [\lambda +
\frac{nH}{2}]\omega_{A2}^{n}
- \mu \omega_{A1}^{n}$\\

\hline
27 & $e_{2}$ & $e_{B}$ & $e_{1}$ & $-e_{B}(\mu)= [\lambda_{n-1}-\lambda]\omega_{2B}^{1}+ \mu \omega^{2}_{2B}$\\

 28 & $e_{2}$ & $e_{B}$ & $e_{2}$ & $
e_{B}(\lambda)=[\lambda_{n-1}-\lambda]\omega_{2B}^{2}- \mu
\omega^{1}_{2B}
$\\

29 & $e_{2}$ & $e_{B}$ & $e_{A}$ & $  (\lambda_{n-1}-\lambda_{3})\omega_{2B}^{A}= [\lambda-\lambda_{3}]\omega_{B2}^{A}- \mu \omega^{A}_{B1}$\\

30 & $e_{2}$ & $e_{B}$ & $e_{B}$ & $ e_{2}(\lambda_{n-1})= [\lambda-\lambda_{n-1}]\omega_{B2}^{B}- \mu \omega^{B}_{B1}$\\

31 & $e_{2}$ & $e_{B}$ & $e_{\widetilde{B}}$ & $[\lambda-\lambda_{n-1}]\omega_{B2}^{\widetilde{B}}- \mu \omega^{\widetilde{B}}_{B1}=0$\\

 32 & $e_{2}$ & $e_{B}$ & $e_{n}$ &
$[\lambda_{n-1}+\frac{nH}{2}]\omega_{2B}^{n} = [\lambda +
\frac{nH}{2}]\omega_{B2}^{n}
- \mu \omega_{B1}^{n}$\\

\hline
33 & $e_{2}$ & $e_{n}$ & $e_{1}$ & $ -(\lambda + \frac{nH}{2})\omega_{2n}^{1}+ \mu \omega_{2n}^{2}= -e_{n}(\mu)$\\

 34 & $e_{2}$ & $e_{n}$ & $e_{2}$ & $-(\lambda +
\frac{nH}{2})\omega_{2n}^{2}- \mu
\omega_{2n}^{1}=e_{n}(\lambda)$\\

35 & $e_{2}$ & $e_{n}$ & $e_{n}$ & $ (\lambda + \frac{nH}{2})\omega_{n2}^{n}- \mu \omega_{n1}^{n}=0$\\

\hline 36. & $e_{A}$ & $e_{B}$ & $e_{1}$ & $(\lambda_{n-1} -
\lambda)\omega_{AB}^{1}+ \mu \omega_{AB}^{2}=
[\lambda_{3}-\lambda]\omega_{BA}^{1}+ \mu \omega_{BA}^{2}$\\

 37 & $e_{A}$ & $e_{B}$ & $e_{2}$ & $(\lambda_{n-1} -
\lambda)\omega_{AB}^{2}- \mu \omega_{AB}^{1}=
[\lambda_{3}-\lambda]\omega_{BA}^{2}- \mu \omega_{BA}^{1}$\\

 38 & $e_{A}$ & $e_{B}$ & $e_{n}$ &
$[\lambda_{n-1}+\frac{nH}{2}]\omega_{AB}^{n}=
[\lambda_{3}+\frac{nH}{2}]\omega_{BA}^{n}$\\

39 & $e_{A}$ & $e_{B}$ & $e_{\widetilde{B}}$ & $ \omega_{BA}^{\widetilde{B}}=0$\\

40 & $e_{A}$ & $e_{B}$ & $e_{\widetilde{A}}$ & $\omega_{AB}^{\widetilde{A}}= 0$\\

\hline 41. & $e_{A}$ & $e_{n}$ & $e_{1}$ & $-(\lambda +
\frac{nH}{2})\omega_{An}^{1}+ \mu \omega_{An}^{2}=
[\lambda_{3}-\lambda]\omega_{nA}^{1}+ \mu \omega_{nA}^{2}$\\

 42 & $e_{A}$ & $e_{n}$ & $e_{2}$ & $-(\lambda +
\frac{nH}{2})\omega_{An}^{2}- \mu \omega_{An}^{1}=
[\lambda_{3}-\lambda]\omega_{nA}^{2}- \mu \omega_{nA}^{1}$\\

43 & $e_{A}$ & $e_{n}$ & $e_{A}$ & $e_{n}(\lambda_{3})= -[\frac{nH}{2}+ \lambda_{3}]\omega_{An}^{A}$\\

 44 & $e_{A}$ & $e_{n}$ & $e_{B}$ &
$-[\lambda_{n-1}+\frac{nH}{2}]\omega_{An}^{B}=
[\lambda_{3}-\lambda_{n-1}]\omega_{nA}^{B}$\\

45 & $e_{A}$ & $e_{n}$ & $e_{\widetilde{A}}$ & $ \omega_{An}^{\widetilde{A}}=0$ \\

46 & $e_{A}$ & $e_{n}$ & $e_{n}$ & $ \omega_{nA}^{n}=0$\\

\hline 47 & $e_{B}$ & $e_{n}$ & $e_{1}$ & $-(\lambda +
\frac{nH}{2})\omega_{Bn}^{1}+ \mu \omega_{Bn}^{2}=
[\lambda_{n-1}-\lambda]\omega_{nB}^{1}+ \mu \omega_{nB}^{2}$\\

 48 & $e_{B}$ & $e_{n}$ & $e_{2}$ & $-(\lambda +
\frac{nH}{2})\omega_{Bn}^{2}- \mu \omega_{Bn}^{1}=
[\lambda_{n-1}-\lambda]\omega_{nB}^{2}- \mu \omega_{nB}^{1}$\\

 49 & $e_{B}$ & $e_{n}$ & $e_{A}$ &
$-[\lambda_{3}+\frac{nH}{2}]\omega_{Bn}^{A}=
[\lambda_{n-1}-\lambda_{3}]\omega_{nB}^{A}$\\

50 & $e_{B}$ & $e_{n}$ & $e_{B}$ & $e_{n}(\lambda_{n-1})= -[\frac{nH}{2}+ \lambda_{n-1}]\omega_{Bn}^{B}$\\

51 & $e_{B}$ & $e_{n}$ & $e_{\widetilde{B}}$ & $ \omega_{Bn}^{\widetilde{B}}=0$ \\

52 & $e_{B}$ & $e_{n}$ & $e_{n}$ & $ \omega_{nB}^{n}=0$\\ [1ex] 
\hline 
\end{longtable}

\medskip

 By using T5, T38, T49, (\ref {3.7}) and (\ref {3.5}) we have
\begin{equation}\label{3.10}
         \omega_{AB}^{n}=\omega_{BA}^{n}=\omega_{An}^{B}=\omega_{AB}^{n}=\omega_{nB}^{A}=\omega_{nA}^{B}=0,
         \quad \omega_{11}^{n}= - \omega_{22}^{n}.
\end{equation}

 Equating T18, T34, T19, T33 and using (\ref {3.5}), we find
  \begin{equation}\label{3.11}
         \omega_{22}^{n}=\omega_{11}^{n}, \quad \omega_{12}^{n}= -\omega_{21}^{n},
\end{equation}
which by use of (\ref {3.7}), (\ref {3.10}) and (\ref {3.5}) give

 \begin{equation}\label{3.12}
         \omega_{22}^{n}=\omega_{11}^{n}= \omega_{12}^{n}= \omega_{21}^{n}=\omega_{2n}^{1}=
         \omega_{1n}^{2}=\omega_{2n}^{2}=
         \omega_{1n}^{1}=0.
\end{equation}

Similarly, using T6, T22, T7, T21, T12, T28, T13, T27 and (\ref
{3.5}), we find
\begin{equation}\label{3.13}
         \omega_{22}^{A}=\omega_{11}^{A}, \quad \omega_{12}^{A}= -\omega_{21}^{A}, \quad \omega_{22}^{B}=\omega_{11}^{B},\quad \omega_{12}^{B}=
         -\omega_{21}^{B}.
\end{equation}

 Using T20, T35, T46, T52 and (\ref {3.5}), we get
\begin{equation}\label{3.14}
         \omega_{n1}^{n}=\omega_{n2}^{n}= \omega_{nn}^{1}=
         \omega_{nn}^{2}= \omega_{nn}^{A}=
         \omega_{nn}^{B}=0.
\end{equation}

 From T39, T40  and (\ref {3.5}), we find
\begin{equation}\label{3.15}
         \omega_{B\widetilde{B}}^{A}=\omega_{A\widetilde{A}}^{B}=0.
\end{equation}

Solving T11, T26, T17, T32  by using (\ref {3.7}) and (\ref {3.5}),
we obtain
\begin{multline}\label{3.16}
         \omega_{1A}^{n}=\omega_{2A}^{n}=\omega_{A1}^{n}=
         \omega_{A2}^{n}= \omega_{1n}^{A}=\omega_{2n}^{A}=\omega_{An}^{1}=
         \omega_{An}^{2}=0, \quad
         \mbox{and}\\
         \omega_{1B}^{n}=\omega_{2B}^{n}=\omega_{B1}^{n}=
         \omega_{B2}^{n}= \omega_{1n}^{B}=\omega_{2n}^{B}=\omega_{Bn}^{1}=
         \omega_{Bn}^{2}=0.\\
\end{multline}

 Using T9, T24, T16, T31, T45, T51  and (\ref {3.5}), we get
\begin{equation}\label{3.17}
         \omega_{A1}^{\widetilde{A}}=\omega_{A2}^{\widetilde{A}}=\omega_{B1}^{\widetilde{B}}=\omega_{B2}^{\widetilde{B}}=\omega_{A\widetilde{A}}^{1}
         =\omega_{A\widetilde{A}}^{2}=\omega_{B\widetilde{B}}^{2}=\omega_{A\widetilde{A}}^{n}=\omega_{B\widetilde{B}}^{n}.
\end{equation}

Now, solving T41, T42, T47 and T48  by using (\ref {3.16}) and (\ref
{3.5}), we obtain
\begin{equation}\label{3.18}
         \omega_{nA}^{1}=\omega_{nA}^{2}=\omega_{n1}^{A}=
         \omega_{n2}^{A}=\omega_{nB}^{1}=\omega_{nB}^{2}=\omega_{n1}^{B}=
         \omega_{n2}^{B}=0.
\end{equation}

Equating T10, T14  by using (\ref {3.5}) and solving with T36, we
get

\begin{equation}\label{3.19}
         (\lambda_{3}-\lambda)\omega_{BA}^{1}=(\lambda_{n-1}-\lambda)\omega_{AB}^{1}, \quad \omega_{AB}^{2}=
         \omega_{BA}^{2}.
\end{equation}

Similarly, equating T25, T29  by using (\ref {3.5}) and solving with
T37 we get

\begin{equation}\label{3.20}
         (\lambda_{3}-\lambda)\omega_{BA}^{2}=(\lambda_{n-1}-\lambda)\omega_{AB}^{2}, \quad \omega_{AB}^{1}=
         \omega_{BA}^{1}.
\end{equation}

Combining (\ref {3.19}) and (\ref {3.20}), and using (\ref {3.5}),
we obtain

\begin{equation}\label{3.21}
         \omega_{BA}^{1}=\omega_{AB}^{1}= \omega_{AB}^{2}=
         \omega_{BA}^{2}=\omega_{B1}^{A}=\omega_{A1}^{B}= \omega_{A2}^{B}=
         \omega_{B2}^{A}=0.
\end{equation}

\vspace{.5 cm}

From the above computations we obtain the following:
\begin{lem}\label{lemma3.1}
 Let $M^{n}_{1}$ be a Lorentz hypersurface in $E^{n+1}_{1}$, having the shape operator (\ref {2.11}) with five distinct eigenvalues
  with respect to a suitable orthonormal basis  $\{e_{1},
e_{2}, \dots , e_{n}\}$. If grad$H$ is space like and in the direction
of $e_{n }$, then

\begin{center}
      $$\nabla_{e_{1}}e_{1}=\sum_{p\neq 1,n}\omega_{11}^{p}e_{p},\hspace{.2 cm}
      \nabla_{e_{1}}e_{2}=\sum_{p\neq 2,n}\omega_{12}^{p}e_{p},\hspace{.2 cm}
      \nabla_{e_{1}}e_{A}=\sum_{p\neq A,n}\omega_{1A}^{p}e_{p},\hspace{.2 cm}
      \nabla_{e_{1}}e_{n}=0,$$\hspace{.2 cm}
      $$\nabla_{e_{2}}e_{1}= \sum_{p\neq 1,n}\omega_{21}^{p}e_{p},\hspace{.2 cm}
      \nabla_{e_{2}}e_{2}=\sum_{p\neq 2,n}\omega_{22}^{p}e_{p},\hspace{.2 cm}
      \nabla_{e_{2}}e_{A}=\sum_{p\neq A,n}\omega_{2A}^{p}e_{p},\hspace{.2 cm}
      \nabla_{e_{2}}e_{n}=0,$$\hspace{.2 cm}
      $$\nabla_{e_{A}}e_{1}=\sum_{p\neq 1,\widetilde{A},B,n}\omega_{A1}^{p}e_{p},\hspace{.2 cm}
      \nabla_{e_{A}}e_{2}=\sum_{p\neq 2,\widetilde{A},B,n}\omega_{A2}^{p}e_{p},\hspace{.2 cm}
      \nabla_{e_{A}}e_{A}=\sum_{p\neq A}\omega_{AA}^{p}e_{p},$$\hspace{.2 cm}
      $$\nabla_{e_{A}}e_{\widetilde{A}}=\sum_{p\neq 1,2,\widetilde{A},B,n}\omega_{A\widetilde{A}}^{p}e_{p},\hspace{.2 cm}
      \nabla_{e_{A}}e_{B}=\sum_{p\neq 1,2,\widetilde{A},B,n}\omega_{AB}^{p}e_{p},\hspace{.2 cm}
      \nabla_{e_{A}}e_{n}= \sum_{p\neq 1,2,\widetilde{A},B,n}\omega_{An}^{p}e_{p},$$\hspace{.2 cm}
      $$\nabla_{e_{B}}e_{1}=\sum_{p\neq 1,A,\widetilde{B},n}\omega_{B1}^{p}e_{p},\hspace{.2 cm}
      \nabla_{e_{B}}e_{2}=\sum_{p\neq 2,A,\widetilde{B},n}\omega_{B2}^{p}e_{p},\hspace{.2 cm}
      \nabla_{e_{B}}e_{A}=\sum_{p\neq 1,2,A,\widetilde{B},n}\omega_{BA}^{p}e_{p},$$\hspace{.2 cm}
      $$\nabla_{e_{B}}e_{\widetilde{B}}=\sum_{p\neq 1,2,A,\widetilde{B},n}\omega_{B\widetilde{B}}^{p}e_{p},\hspace{.2 cm}
      \nabla_{e_{B}}e_{B}=\sum_{p\neq B}\omega_{BB}^{p}e_{p},\hspace{.2 cm}
      \nabla_{e_{B}}e_{n}= \sum_{p\neq 1,2,A,\widetilde{B},n}\omega_{Bn}^{p}e_{p},$$\hspace{.2 cm}
      $$\nabla_{e_{n}}e_{1}=\omega_{n1}^{2}e_{2},\hspace{.2 cm}
      \nabla_{e_{n}}e_{2}=\omega_{n2}^{1}e_{1},\hspace{.2 cm}
      \nabla_{e_{n}}e_{A}=\sum_{p\neq 1,2,A,B,n}\omega_{nA}^{p}e_{p},$$\hspace{.2 cm}
      $$\nabla_{e_{n}}e_{B}=\sum_{p\neq 1,2,A,B,n}\omega_{nB}^{p}e_{p},\hspace{.2 cm}
      \nabla_{e_{n}}e_{n}=0,\hspace{.2 cm}
      \nabla_{e_{2}}e_{B}=\sum_{p\neq B,n}\omega_{2B}^{p}e_{p},$$\hspace{.2 cm}
\end{center}
\end{lem}

\section{\textbf{Proof of Theorem \ref{main}}}
We use Lemma \ref {lemma3.1}, (\ref {2.5}) and (\ref {3.5}) to
evaluate $g(R(e_{2},e_{n})e_{2},e_{n})$,  and we obtain
\begin{equation}\label{4.1}
\frac{nH}{2} \lambda =0, \quad\mbox{so}\quad \lambda=0,
\end{equation}
as $H\neq0$.

Now, using (\ref {2.4}), we have trace$S=nH$ and using (\ref {4.1})
and the value of $\lambda_{n}=-\frac{nH}{2}$, we get

\begin{equation}\label{4.2}
(r-2)\lambda_{3} +(n-r-1)\lambda_{n-1}=\frac{3nH}{2}.
\end{equation}

Using T43,  T50, (\ref {4.2}) and (\ref {3.5}), we obtain

\begin{equation}\label{4.3}
3ne_{n}(H)=[n(n-r+2)H-2(r-2)\lambda_{3}]\omega^{n}_{BB}
+(r-2)(2\lambda_{3}+nH)\omega^{n}_{AA}.
\end{equation}

Using Lemma \ref {lemma3.1}, (\ref {3.6}) and the fact that
$[e_{i}\hspace{.1 cm}
e_{n}](H)=0=\nabla_{e_{i}}e_{n}(H)-\nabla_{e_{n}}e_{i}(H),$ for $i
=1, 2, \ldots, n-1$, we obtain
      \begin{equation}\label{4.4}
      e_{i}e_{n}(H)= 0.
      \end{equation}

Using T3, T6, T4, T12, (\ref {3.13}), (\ref {4.1}) and (\ref {3.5}),
we find that

\begin{center}
         $\omega_{12}^{A}(\lambda_{3}^{2}-
         \mu^{2})=0 \quad \mbox{and} \quad
         \omega_{12}^{B}(\lambda_{n-1}^{2}-
         \mu^{2})=0.$
\end{center}

\vspace{.5 cm}

Therefore, we need to consider the following cases:\\

\noindent
\textbf{Case A.} \ $\lambda_{3}^{2}= \mu^{2}$, $ \lambda_{n-1}^{2}=
\mu^{2}$.\\
In this case we have that
$e_{n}(\lambda_{3})=e_{n}(\lambda_{n-1})=0$, because from T33, it is
$e_{n}(\mu)=0$. Using T43 and T50, we obtain that
$\omega^{n}_{AA}=\omega^{n}_{BB}=0$. Therefore, from (\ref {4.3}),
we have
$e_{n}(H)=0$, which contradicts that $e_{n}(H)\ne0.$\\

\noindent
\textbf{Case B.} \ $ \omega^{A}_{12}=0$, $ \lambda_{n-1}^{2}=
\mu^{2}$.\\
In this case we have that $e_{n}(\lambda_{n-1})=0$, which implies
from T50  that $\omega^{n}_{BB}=0$. Therefore, using (\ref {2.5})
and Lemma \ref {lemma3.1} to evaluate
$g(R(e_{B},e_{n})e_{n},e_{B})$, we get
      \begin{equation}\label{4.5}
      -\frac{nH}{2}\lambda_{n-1}=0, \quad \mbox{so} \quad \lambda_{n-1}=0.
      \end{equation}

 Using (\ref {4.2}), (\ref {4.5}) and (\ref {4.3}), we obtain

\begin{equation}\label{4.6}
3e_{n}(H)=(r+1)H\omega^{n}_{AA},
\end{equation}
which implies
\begin{equation}\label{4.7}
       \omega_{33}^{n}=\omega_{44}^{n}=\dots=\omega_{rr}^{n}.
      \end{equation}

 Also, using that trace$S^{2}=(r-2)\lambda_{3}^{2}+(n-r-1)\lambda_{n-1}^{2}- 2 \mu^2$, (\ref {2.10}), (\ref {4.2}), (\ref {4.5}), (\ref {4.7}) and
 Lemma \ref {lemma3.1}, then equation (\ref {2.8}) with respect to the basis $\{e_{1},
e_{2}, \dots , e_{n}\}$ reduces to
        \begin{equation}\label{4.8}
      -e_{n}e_{n}(H)+(r-2)\omega_{AA}^{n}e_{n}(H)+ H[\frac{n^{2}(r+7)H^{2}}{4(r-2)}- 2\mu^2]=\alpha
      H.
      \end{equation}

      Using (\ref {2.5}), (\ref {3.5}) and Lemma \ref {lemma3.1} to evaluate  $g(R(e_{A},e_{n})e_{n},e_{A})$, we
      get
      \begin{equation}\label{4.9}
      e_{n}(\omega_{AA}^{n})-(\omega_{AA}^{n})^{2}=
      -\frac{3n^{2}H^{2}}{4(r-2)}.
      \end{equation}

        Differentiating (\ref {4.6}) along $e_{n}$ and using (\ref {4.9}) and (\ref {4.6}), we
        get
        \begin{equation}\label{4.10}
        3e_{n}e_{n}(H)=\frac{(r+1)(r+4)H}{3}(\omega^{n}_{AA})^{2}-\frac{3n^{2}(r+1)H^{3}}{4(r-2)}.
        \end{equation}

        Eliminating $e_{n}e_{n}(H)$ and $e_{n}(H)$ from (\ref {4.8}) by using (\ref {4.10}) and (\ref {4.6}), we obtain
        \begin{equation}\label{4.11}
       \frac{2(r+1)(r-5)}{3}(\omega_{AA}^{n})^{2}+\frac{3n^{2}(r+4)H^{2}}{2(r-2)}-6
       \mu^2=3\alpha.
      \end{equation}

       Differentiating (\ref {4.11}) with respect to $e_{n}$ and using (\ref {4.9}) and (\ref {4.6}), we
       find
       \begin{equation}\label{4.12}
      \frac{2(r-5)}{3} (\omega_{AA}^{n})^{2}+\frac{9n^{2}H^{2}}{2(r-2)}=
      0.
      \end{equation}

       Again, acting along $e_{n}$ on (\ref {4.12}) and using (\ref {4.9}) and (\ref {4.6}), we get
        \begin{equation}\label{4.13}
      \frac{4(r-5)}{3} (\omega_{AA}^{n})^{2}+\frac{2n^{2}(r+4)H^{2}}{r-2}=
      0.
      \end{equation}

       Hence, from (\ref {4.12}) and (\ref {4.13}), we obtain that $H$ must be
       zero.

      \smallskip
\noindent
       $\textbf{Case C:}$ \ $\lambda_{3}^{2}= \mu^{2}$,
       $\omega^{B}_{12}=0$.\ Analogously to Case B it can be shown that
$H=0$.

\smallskip
\noindent
 $\textbf{Case D:}$ \ $ \omega^{A}_{12}=0$, $ \omega^{B}_{12}=0$.\\

    Using T3, T4  and (\ref {3.13}), we find that

   \begin{equation}\label{4.14}
         \omega_{21}^{B}=\omega_{21}^{A}=\omega_{22}^{B}=\omega_{22}^{A}=\omega_{11}^{B}=\omega_{11}^{A}=0.
       \end{equation}

       Now, by computing  $g(R(e_{A},e_{1})e_{A},e_{n})$,
$g(R(e_{B},e_{1})e_{B},e_{n})$, $g(R(e_{A},e_{2})e_{A},e_{n})$,
$g(R(e_{B},e_{2})e_{B},e_{n})$ using (\ref {2.5}), (\ref {3.5}) and
Lemma \ref {lemma3.1}, we obtain that

       \begin{equation}\label{4.15}
         e_{1}(\omega_{AA}^{n})-\omega_{AA}^{n}\omega_{AA}^{1}=0,
       \end{equation}

       \begin{equation}\label{4.16}
         e_{1}(\omega_{BB}^{n})-\omega_{BB}^{n}\omega_{AA}^{1}=0,
       \end{equation}

\begin{equation}\label{4.17}
         e_{2}(\omega_{AA}^{n})-\omega_{AA}^{n}\omega_{AA}^{2}=0,
       \end{equation}

\begin{equation}\label{4.18}
         e_{2}(\omega_{BB}^{n})-\omega_{BB}^{n}\omega_{BB}^{2}=0,
       \end{equation}
       respectively.

       Differentiating (\ref {4.2}) along $e_{1}$ and using T8, T15  and (\ref {4.2}), we obtain

\begin{equation}\label{4.19}
         2(r-2)[\lambda_{3}\omega_{AA}^{1}-\mu\omega_{AA}^{2}]+[(3nH-2(r-2)\lambda_{3})\omega_{BB}^{1}-2\mu(n-r-1)\omega_{BB}^{2}]=0.
       \end{equation}

      Similarly, differentiating (\ref {4.2}) along $e_{2}$ and using T23, T30 and (\ref {4.2}), we obtain that

\begin{equation}\label{4.20}
         2(r-2)[\lambda_{3}\omega_{AA}^{2}+\mu\omega_{AA}^{1}]+[(3nH-2(r-2)\lambda_{3})\omega_{BB}^{2}+2\mu(n-r-1)\omega_{BB}^{1}]=0.
       \end{equation}

       Multiplying (\ref {4.19}) and (\ref {4.20}) by $\lambda_{3}$ and $\mu$
       respectively, and then adding, we get

       \begin{multline}\label{4.21}
         2(r-2)(\lambda_{3}^{2}+\mu^{2})\omega_{AA}^{1}+[\{\lambda_{3}(3nH-2(r-2)\lambda_{3})+2\mu^{2}(n-r-1)\}\omega_{BB}^{1}\\
         +\{\mu(3nH-2(r-2)\lambda_{3})-2\mu \lambda_{3}(n-r-1)\}\omega_{BB}^{2}]=0.
       \end{multline}

        Now, multiplying (\ref {4.19}) and (\ref {4.20}) by $\mu$ and $\lambda_{3}$
       respectively, and subtracting, we get

       \begin{multline}\label{4.22}
         2(r-2)(\lambda_{3}^{2}+\mu^{2})\omega_{AA}^{2}+[\{\lambda_{3}(3nH-2(r-2)\lambda_{3})+2\mu^{2}(n-r-1)\}\omega_{BB}^{2}\\
         +\{-\mu(3nH-2(r-2)\lambda_{3})+2\mu \lambda_{3}(n-r-1)\}\omega_{BB}^{1}]=0.
       \end{multline}

       Differentiating (\ref {4.3}) along $e_{1}$ and using (\ref {4.4}), (\ref {4.15}),
       (\ref {4.16}), we obtain

       \begin{multline}\label{4.23}
         2(r-2)[\lambda_{3}\omega_{AA}^{1}-\mu\omega_{AA}^{2}](\omega^{n}_{AA}-\omega^{n}_{BB})+(n(n-r+2)H-2(r-2)\lambda_{3})\omega_{BB}^{n}\omega_{BB}^{1}\\
         +(r-2)(2\lambda_{3}+nH)\omega_{AA}^{n}\omega_{AA}^{1}=0.
       \end{multline}

      Similarly, differentiating (\ref {4.3}) along $e_{2}$ and using (\ref {4.4}), (\ref {4.17}),
       (\ref {4.18}), we obtain

       \begin{multline}\label{4.24}
         2(r-2)[\lambda_{3}\omega_{AA}^{2}+\mu\omega_{AA}^{1}](\omega^{n}_{AA}-\omega^{n}_{BB})+(n(n-r+2)H-2(r-2)\lambda_{3})\omega_{BB}^{n}\omega_{BB}^{2}\\
         +(r-2)(2\lambda_{3}+ nH)\omega_{AA}^{n}\omega_{AA}^{2}=0.
       \end{multline}

       Eliminating $\omega_{AA}^{1}$ and $\omega_{AA}^{2}$ from
       (\ref {4.23}) using (\ref {4.19}) and (\ref {4.21}), we obtain

       \begin{multline}\label{4.25}
        \omega_{AA}^{n}[2\mu(n-r-1)\omega^{2}_{BB}-(3nH-2(r-2)\lambda_{3})\omega_{BB}^{1}-\frac{(2\lambda_{3}+nH)}{2(\lambda_{3}^{2}+\mu^{2})}
        (\{\lambda_{3}(3nH-2(r-2)\lambda_{3})+\\2\mu^{2}(n-r-1)\}\omega_{BB}^{1}
         +\{\mu(3nH-2(r-2)\lambda_{3})-2\mu \lambda_{3}(n-r-1)\}\omega_{BB}^{2})]
         +\omega_{BB}^{n}[-2\mu(n-r-1)\omega^{2}_{BB}\\+(n(n-r+5)H-4(r-2)\lambda_{3})\omega_{BB}^{1}]=0.
       \end{multline}

       Similarly, eliminating $\omega_{AA}^{1}$ and $\omega_{AA}^{2}$ from
       (\ref {4.24}) using (\ref {4.20}) and (\ref {4.22}), we obtain

\begin{multline}\label{4.26}
        \omega_{AA}^{n}[-2\mu(n-r-1)\omega^{1}_{BB}-(3nH-2(r-2)\lambda_{3})\omega_{BB}^{2}-\frac{(2\lambda_{3}+nH)}{2(\lambda_{3}^{2}+\mu^{2})}
        \huge{ ( }\{\lambda_{3}(3nH-2(r-2)\lambda_{3})+\\2\mu^{2}(n-r-1)\}\omega_{BB}^{2}
         +\{-\mu(3nH-2(r-2)\lambda_{3})+2\mu \lambda_{3}(n-r-1)\}\omega_{BB}^{1})]
         +\omega_{BB}^{n}[2\mu(n-r-1)\omega^{1}_{BB}\\+(n(n-r+5)H-4(r-2)\lambda_{3})\omega_{BB}^{2}]=0.
       \end{multline}

       Now, eliminating $\omega^{n}_{AA}$ and $\omega^{n}_{BB}$ from
       (\ref {4.25}) and (\ref {4.26}), we get

       \begin{multline}\label{4.27}
        [(\omega_{BB}^{1})^{2}+(\omega_{BB}^{2})^{2}][2PQ(\lambda_{3}^{2}+\mu^{2})+Q(2\lambda_{3}+nH)(\lambda_{3}P-\mu
        R)-2PR(\lambda_{3}^{2}+\mu^{2})\\
        -P(2\lambda_{3}+nH)(\lambda_{3}R+ \mu P)]=0,
       \end{multline}
       where $P=2\mu(n-r-1)$, $Q=n(n-r+5)H-4(r-2)\lambda_{3}$ and
       $R=3nH-2(r-2)\lambda_{3}.$\\

      We now claim that $(\omega_{BB}^{1})^{2}+(\omega_{BB}^{2})^{2}\neq
       0$.\\

      Indeed, if $(\omega_{BB}^{1})^{2}+(\omega_{BB}^{2})^{2}=0$, we have,
      $\omega_{BB}^{1}=\omega_{BB}^{2}=0$ as connection coefficients are real numbers. Then, using (\ref {4.21}) and (\ref {4.22}),
      we have $\omega_{AA}^{1}=\omega_{AA}^{2}=0$.

     Therefore, by computing  $g(R(e_{B},e_{2})e_{B},e_{1})$,
$g(R(e_{A},e_{2})e_{A},e_{1})$, using (\ref {2.5}), (\ref {4.14}),
(\ref {3.5}) and Lemma \ref {lemma3.1}, we obtain

      \begin{equation}\label{4.28}
         \lambda_{n-1}\mu=0,
       \end{equation}

       \begin{equation}\label{4.29}
         \lambda_{3}\mu=0,
       \end{equation}
       respectively, which implies $\lambda_{3}=\lambda_{n-1}=0$.
       Using T43 and T50, we obtain that
       $\omega_{AA}^{n}=\omega_{BB}^{n}=0$. Also, from (\ref {4.3}) we have $e_{n}(H)=0$, which is a contradiction. Hence the claim is proved.\\

Therefore, from (\ref {4.27}) we have that
\begin{multline}\label{4.30}
        f(\lambda_{3},H)\equiv2PQ(\lambda_{3}^{2}+\mu^{2})+Q(2\lambda_{3}+nH)(\lambda_{3}P-\mu
        R)-2PR(\lambda_{3}^{2}+\mu^{2})
        -P(2\lambda_{3}+nH)(\lambda_{3}R+ \mu P)=0.
       \end{multline}

       Now, using T7, T13, T33, (\ref {3.12}) and (\ref {4.14}), we
       obtain

      \begin{equation}\label{4.31}
         e_{A}(\mu)=e_{B}(\mu)=e_{n}(\mu)=0.
       \end{equation}

       Also, from T1, T2  and (\ref {4.1}),
       we have $e_{1}(\mu)=e_{2}(\mu)=0$, which implies from (\ref {4.31})
       that $\mu$ is constant in each direction.

       Differentiating (\ref {4.30}) along $e_{1}$ and $e_{2}$, we have

       \begin{equation}\label{4.32}
         e_{1}(\lambda_{3})(g(\lambda_{3},H))=0,
       \end{equation}
       and
\begin{equation}\label{4.33}
         e_{2}(\lambda_{3})(g(\lambda_{3},H))=0,
       \end{equation}
       respectively, where $g(\lambda_{3},H)= 4P\lambda_{3}(Q-R)-4P(\lambda_{3}^{2}+\mu^{2})(r-2)+2(PQ\lambda_{3}-QR\mu-\lambda_{3}PR-P^{2}\mu)
       +(2\lambda_{3}+nH)(PQ-2\lambda_{3}P(r-2)+2(r-2)(2R+Q)\mu-PR)$.

Now, if $g(\lambda_{3},H)\neq 0$, we have $e_{1}(\lambda_{3})=0$ and
$e_{2}(\lambda_{3})=0$ which implies from T8, T15, T23 and T30 that
$\omega_{BB}^{1}=\omega_{BB}^{2}=\omega_{AA}^{1}=\omega_{AA}^{2}=0$.
As we have already proved from (\ref {4.28}) and (\ref {4.29}) this
gives a contradiction.

Therefore, we have

\begin{equation}\label{4.34}
         g(\lambda_{3},H)= 0,
       \end{equation}
  which is a polynomial equation in $\lambda_{3}$ and $H$.

We rewrite $f(\lambda_{3},H)$, $g(\lambda_{3},H)$ as polynomials
$f_{H}(\lambda_{3}), g_{H}(\lambda_{3})$ of $\lambda_{3}$ with
coefficients in the polynomial ring $R[H]$ over $\mathbb{R}$. Since
$f_{H}(\lambda_{3}) = g_{H}(\lambda_{3}) = 0$,  $\lambda _3$ is a
common root of $f_H, g_H$, hence by Lemma \ref{lemma2.1} it is
$\Re(f_{H}, g_{H}) = 0$. It is obvious that $\Re(f_{H}, g_{H})$ is a
polynomial of $H$ with constant coefficients, therefore $H$ must be
a constant. This contradicts the first of relations (\ref {3.6}).

Cases A, B, C, D conclude the proof of Theorem \ref {main}. \qed

\medskip
\noindent
\textbf{Acknowledgement}. \small \emph{ The first author is
grateful to Guru Gobind Singh Indraprastha University for providing
IPRF fellowship  to pursue research. The second author was supported
by Grant $\# E.037$ from the Research Committee of the University of Patras  (Programme K. Karatheodori).}



Author's address:\\
\\
\textbf{Deepika}\\
University School of Basic and Applied Sciences,\\
Guru Gobind Singh Indraprastha University,\\
Sector-16C, Dwarka, New Delhi-110078, India.\\
\textbf{Email:} sdeep2007@gmail.com\\
\\
\textbf{Andreas Arvanitoyeorgoes}\\
University of Patras,\\
Department of Mathematics,\\
 GR-26500 Patras, Greece.\\
 \textbf{Email:} arvanito@math.upatras.gr\\
\\
\textbf{Ram Shankar Gupta}\\
University School of Basic and Applied Sciences,\\
 Guru Gobind Singh Indraprastha University,\\
Sector-16C, Dwarka, New Delhi-110078, India.\\
\textbf{Email:} ramshankar.gupta@gmail.com

\end{document}